\documentstyle[amssymb,12pt]{article}  
\title{Super-connections and non-commutative geometry}
\author{ Victor Nistor\thanks{Partially supported by the National Science
Foundation under grant number 
DMS-9205542, an NSF Young Investigator Award
DMS-9457859 and a Sloan Research Fellowship, manuscripts available
from {\bf 
http:{\scriptsize//}www.math.psu.edu{\scriptsize/}nistor{\scriptsize/}}}
\\
{\normalsize Pennsylvania State University} \vspace{-0.3cm}\\
{\normalsize Department of Mathematics} \vspace{-0.3cm}\\
{\normalsize McAllister Bldg.} \vspace{-0.3cm}\\
{\normalsize University Park, PA 16802}}
\begin{document}
\input{amssym.def}
\maketitle
\begin{abstract}
{We show that Quillen's formalism for computing
the Chern character of the index using superconnections extends 
to arbitrary operators with functional calculus.
We thus remove the condition that the operators have, up to homotopy,
a gap in the spectrum. This is proved  using differential graded algebras
and non-commutative differential forms.
Our results also give a new proof of the coincidence of 
the Chern character of a difference bundle
defined using super-connections with the classical definition.}
\end{abstract}

\newtheorem{theorem}{Theorem.}[section]
\newtheorem{maintheorem}[theorem]{Main Theorem}
\newtheorem{lemma}[theorem]{Lemma.}
\newtheorem{proposition}[theorem]{Proposition.}
\newtheorem{corollary}[theorem]{Corollary.}
\newtheorem{definition}[theorem]{Definition.}
\newtheorem{mainlemma}[theorem]{Main Lemma}

\newenvironment{proof}{\begin{trivlist}{\item[] \bf Proof.\/}}{\end{trivlist}}
\newenvironment{example}{\begin{trivlist}{\item[] \bf Examples\/}}{\end{trivlist}}

\newcommand{\Om}{\Omega^*}
\newcommand{\Sy}{\Sigma^\infty}
\newcommand{\Di}{{\cal D}}
\newcommand{\Tt}{\mbox{$a_{0}\otimes \ldots \otimes a_n$}}
\newcommand{\Dir}{{/\hspace{-3mm}D}}
\newcommand{\supp}{\mbox{\rm supp}}
\newcommand{\D}{\;/\hspace{-3mm}D}
\newcommand{\sqep}{$\;\;\sqcap \!\!\!\! \sqcup$}
\newcommand{\cycle}{(\Omega, \nabla, \tau)}
\newcommand{\cycletil}{(\tilde{\Omega}, d, \tilde{\tau})}

\tableofcontents

\section{Introduction}
This paper is a sequel of \cite{NiD} where the well known
McKean-Singer formula was generalized. 
The problem we are concerned with is to
find explicit formulae for the Connes-Karoubi character of 
the index of an `unbounded Fredholm family $D$
parametrized by an algebra $B$'. The usual family index theorem corresponds
to the case $B=C^\infty(Y)$ where $Y$ is the parameter space.

Given such an unbounded Fredholm family  $D$ of hermitian
operators on a Hilbert space, its 
$K$-theoretical index is 
defined using the graph projection of $D$:
\begin{eqnarray}
\label{eqgraph}
p=\left[\begin{array}{cc}1- e^{-P^*P}&\tau(P^*P)P^*\\ 
\tau(PP^*)P & e^{-PP^*}
\end{array}\right]
\end{eqnarray}
where 
\begin{equation}
\label{eqD}
D=\left[\begin{array}{cc}0 & P^*\\ P & 0
\end{array}\right]
\end{equation} 
and $\tau$ is a smooth  even  function satisfying
$\tau(x)^2x^2=e^{-x^2}(1-e^{-x^2})$. 
This projection is also called the Bott or the Wasserman projection
by some authors.

In the case of a single Fredholm 
operator it is easy to see that the relative dimension 
of the graph projection and the constant projection $e_0$
$$
e_0=\left[\begin{array}{cc}1 & 0\\ 0 & 0
\end{array}\right]
$$ 
is the same as the index of $P$.
For families the difference of the K-theory classes $[p]-[e_0]\in 
K_0(C^\infty(Y))=K^0(Y)$ coincides with the the index bundle
defined by Atiyah-Singer \cite{AS4}.

This provides us with a method of computing 
the Connes-Karoubi classes of the index
in cyclic homology. 
Computations using the graph projections were carried out
in \cite{CM,Nest-Tsygan}.
Quite often however it is easier to work with heat kernels rather than 
the above projection
as in \cite{ABP,Bismut,Getzler,Lott}.
For the families index theorem
\cite{AS4} this requires superconnections \cite{Quillen0}. 
As for connections, to a superconnection there is associated
the supercurvature endomorphism which is a closed form, 
nonhomogeneous in general.
The main fact is then that the Chern character of the index 
coincides with (the cohomology class of) the exponential of the supercurvature
$$
ch(\mbox{Ind}(D))=exp(-(tD+\nabla)^2). 
$$
The proof of this in the case of families parametrized by
a {\em compact} manifold, goes as follows. 
One first establishes this principle for operators with gaps 
in the spectrum and the same index bundle. 
Then, using the invariance under homotopy 
of the cohomology class of the exponential of the supercurvature,
the general case is reduced to the above mentioned one. 
This procedure
is due to Bismut \cite{Bismut} and was also used in \cite{Lott}.

There are situations, however, where there are no operators 
with gaps in the 
spectrum, but superconnections and the exponential of the 
supercurvature still make sense. A natural problem is what they
represent. Situation of this sort appear in the study of 
foliations, see \cite{Heitsch}, for families of b-pseudodifferential
operators \cite{BiC,MP2}, or on open manifolds.

The main result of this paper establishes the equality of
the Connes-Karoubi character of the graph projection and of the
Chern character of the superconnection whenever some sort of 
pseudodifferential calculus exists. 
See Theorem \ref{thmMain}. The assumptions of our 
theorem are verified in the classical cases, for pseudodifferential
operators along the leaves of a foliation and for
open manifolds. It also applies to 
algebraic (or formal) settings, such as the case of formal 
deformations $B=C_c^\infty(Y)[[h]]$ of certain
commutative algebras.

The main ingredients of the proof are to construct a certain 
completion of the universal differential algebra of
$C_c^\infty(\Bbb R)$ and to prove the theorem in this
formal but universal setting. The proof relies heavily
on cyclic cohomology. It is interesting to mention
that there are classical cases when our results can be formulated
without any reference to noncommutative geometry, but for 
which no classical proof is available. 

Let us mention the following example. Let $X$ be an open 
$\sigma$-compact manifold.
The $K$-theory groups of $X$ with compact support can be 
defined in two 
equivalent ways, either as the kernel of $K^*(X\cup\{\infty\})\to 
K^*(\{\infty\})$, or as equivalence classes of triples
$(E_0,E_1,P)$ where $P\in \mbox{End}(E_0,E_1)$ is invertible outside 
a compact set. We can think of  the difference 
bundle $[E_0]-[E_1] \in K^0(X)$ as ``the index of $P$''.
We can assume that $E_0$ and $E_1$ are endowed with hermitian
metrics and that $||P||\to \infty$ at $\infty$. To such a 
triple $(E_0,E_1,P)$ there is associated an element 
$[p]-[e_0]\in \mbox{ker}(K^*(X\cup\{\infty\})\to  K^*(\{\infty\}))$
where $p$ and $e_0$ are the projections defined above
(see definition \ref{def4.2}. Our main result
implies that the super-connection Chern character as defined
by Quillen \cite{Quillen0} coincides with the classical 
Chern character $K^*(X)\subset K^*(X\cup \{\infty\}) \to 
H^*(X\cup \{\infty\})\otimes \Bbb R$. This result, implicit in 
Quillen's work, was first published by Berline and Vergne \cite{Berline-Vergne}. 

Other examples and applications are to the Chern character of Dirac
operators on foliated manifolds, leading to the vanishing of certain
secondary characteristic classes for foliations with positive scalar
curvature along the leaves, see \cite{Heitsch,NiI}.

Our main result extends to the odd case, the proof being the
same, word for word. This gives a method to compute the Connes-Karoubi
character on $K_1$ using super-connections. In case $n=1$  this
recovers the main result from \cite{Getzler-top}. Some applications
of this theorem to $1$-dimensional foliations will be included elsewhere.

Some very interesting results related to the results in this paper 
are contained in \cite{CM3} where $Diff$-invariant structures 
and a `universal' index theorem are
treated in detail. See also \cite{CoR}.

I would like to thank Mich\`{e}le Vergne for bringing her joint work 
with Nicole Berline \cite{Berline-Vergne} to my attention.

\section{Connections and traces}
We  consider following Connes \cite{Co2}
the construction of cyclic cocycles from algebras
endowed with connections and traces
that vanish on covariant derivatives
(i.e. from what are called bellow `cycles'). 

A similar
problem was solved by Quillen in the framework 
of the $(b,B)$-bicomplex \cite{Quillen2}. We review and extend 
here Connes's approach to  the superalgebra case,
carefully keeping track of signs. 
The terms {\em graded vector space}, $\Bbb
Z_2$-{\em graded vector space} and {\em super vector space} 
will be used interchangeably. 

Unless otherwise mentioned all the spaces considered will be 
locally convex topological
spaces, all algebras will be  locally convex algebras with jointly continuous
multiplications, and all linear maps will be assumed continuous.

We begin by considering a filtered superalgebra
$\Omega=F_0\Omega\supset F_{-1}\Omega
\supset \ldots $ endowed with an odd  graded derivation $\nabla$
satisfying $\nabla(F_{-k}\Omega)\subset
F_{-k-1}\Omega$. 
The $\Bbb Z_2$-degree on an element $a\in \Omega$ will be denoted
by $\partial a$.
Thus we have that $F_{-k}\Omega$ is $\Bbb Z_2$-graded vector space  and that
$\nabla (ab)=(\nabla a)b+(-1)^{\partial a}a\nabla b$.
The multiplication 
$F_{-k}\Omega\otimes F_{-l}\Omega \to F_{-k-l}\Omega$ is even.
We do not assume $\Omega$ to have a unit.

Recall that a multiplier $(l,r)$ 
of $\Omega$ consists of a pair of 
linear maps $l,r:\Omega \to \Omega$ satisfying $l(ab)=l(a)b$,
$r(ab)=ar(b)$ and $al(b)=r(a)b$ for any $a,b \in \Omega$.
We shall call $\nabla$ a {\em connection} if there exists 
a multiplier $\omega=(l,r)$, called the {\em curvature},
such that $\nabla^2(a)=l(a)-r(a)$
and $[\nabla,l]=[\nabla,r]=0$. If $\Omega$ has a unit multipliers are
in one-to-one correspondence with elements of $\Omega$, so 
$\omega \in F_{-2}\Omega$, $l(a)=\omega a$ and $r(a)=a\omega$. If moreover
the filtration of $\Omega$ comes from a grading
(i.e. $\Omega=\oplus _{k=0}^\infty \Omega^k$) and 
$\nabla(\Omega^k)\subset \Omega^{k+1}$
we have $\omega\in \Omega^2$ and 
the above definition coincides 
with the usual one. In the unital case we also  
have $\nabla^2=ad_\omega$ and 
$\nabla(\omega)=0$, see \cite{Co2}. 
By a {\em supertrace} on $\Omega$ we shall mean
a graded trace with respect to the $\Bbb Z_2$-degree, that
is a map $\tau:\Omega\to \Bbb C$ satisfying 
$\tau(ab)=(-1)^{\partial a \partial b}\tau(ba)$.
A supertrace $\tau$ will be called {\em closed} 
if $\tau(\nabla a)=0$ for any $a\in \Omega$.

\begin{definition}[Connes] 
\label{def1}
 (i) A {\em chain} $(\Omega,\nabla,\tau)$ 
is a filtered superalgebra $\Omega=F_0\Omega\supset F_{-1}\Omega
\supset \ldots $ 
with an odd  connection $\nabla$ satisfying $\nabla(F_{-k}\Omega)\subset
F_{-k-1}\Omega$ and a 
supertrace $\tau:\Omega\to \Bbb C$.

(ii) A  chain $(\Omega,\nabla,\tau)$ will be called {\em exact} if 
$\nabla ^2=0$.

(iii) A chain for which the trace is closed (i.e. $\tau(\nabla a)=0\;\forall \;a$)  
is called a {\em cycle}.
\end{definition}

Recall also \cite{Co2} that the universal differential (bi)graded 
algebra associated to a $\Bbb Z_2$-graded algebra $A$ can be realized
as $\Omega(A)=\oplus_{n=0}^\infty A^+\otimes A^{\otimes n}$
where $A^+=A\oplus \Bbb C$ 
is the algebra with an  adjoint unit (even if $A$ already
had one) and $d(\Tt)=1\otimes \Tt$.
The filtration is 
$F_{-k}\Omega(A)=\oplus_{n=k}^\infty A^+\otimes A^{\otimes n}$. 
See also \cite{Arveson,Kastler}.
Note that the $\Bbb Z_2$-degree of a form in $\Omega^n(A)$ may be different
from $n  \!\!\!\!\pmod {\Bbb Z_2}$ if $A$ is not trivially graded.

Recall \cite{Co2} that a cyclic cocycle $\varphi$ of order $n$ 
on a superalgebra $A$ is a multilinear map 
$\varphi:A^{\otimes n+1}\to \Bbb C$ satisfying 
\begin{eqnarray}
\label{eq2}
\varphi(a_0a_1,\ldots,a_n)-\varphi(a_0,a_1a_2,\ldots,a_n)+
\ldots & + \!\!& (-1)^{n-1}\varphi(a_0,\ldots ,a_{n-1}a_n)
\nonumber \\
+(-1)^{\nu+n}\varphi(a_na_0,a_1,\ldots,a_{n-1})&=&0
\end{eqnarray}
and
\begin{equation}
\label{eq3}
\varphi(a_n,a_0,\ldots,a_{n-1})=(-1)^{\nu+n}\varphi(a_0,a_1,\ldots,a_n)
\end{equation}
where $\nu={\partial a_n(\partial a_0+\ldots
+\partial a_{n-1})}$.

The first equation can also be written as $\varphi\circ b=0$ where
$b$ is the Hochschild boundary \cite{Kassel3}, and hence (\ref{eq2}) 
means that $\varphi$ is a Hochschild cocycle.
The second equation is the cyclic property in the graded case
\cite{Co2,Kassel3}. One can see from the definition that a $0$-cyclic
cocycle on $A$ is simply a supertrace. 
Our definition is taken from
\cite{Kassel3} and 
is different from the one used in \cite{Kastler}.

We record for convenience some easy generalizations 
to the graded case of some  computations in \cite{Co2}. 
\begin{lemma}
\label{lemma1}
Denote by $\Phi(\alpha_0\otimes \alpha_1\otimes \ldots \otimes \alpha_n)=
(-1)^\mu \alpha_0d\alpha_1\ldots d\alpha_n$, then
\begin{equation}
\Phi(b(\alpha_0\otimes \alpha_1\otimes \ldots \otimes \alpha_n))
=(-1)^{n-1}[\Phi(\alpha_0\otimes 
\ldots \otimes \alpha_{n-1}),\alpha_n]
\end{equation}
Here $d^2=0$, $\mu=\sum_{j=0}^n (n-j)\partial a_j=\partial \alpha_{n-1}+\partial 
\alpha_{n-3}+\partial \alpha_{n-5}+\ldots$,
$b$ is the Hochschild boundary and the commutator is the graded commutator:
$[a,a']=aa'-(-1)^{\partial a \partial a'}a'a$.
\end{lemma}
\begin{proof}
Direct computation:
\begin{eqnarray*}
\lefteqn{\Phi(b(\alpha_0\otimes \alpha_1
\otimes \ldots \otimes \alpha_n))=}\\
&\!=\!&(-1)^{(n-1)\partial \alpha_0 + \sum_{j=1}^{n}(n-j)\partial \alpha_j }
\alpha_0\alpha_1d\alpha_2\ldots d\alpha_n\\
&\!+\!&\!\sum_{k=1}^{n-1}
(-1)^{\sum_{0}^{k}(n-j-1)\partial \alpha_j
+ \sum_{k+1}^{n}(n-j)\partial \alpha_j + k+\partial \alpha_k}
\alpha_0d\alpha_1\ldots d\alpha_{k-1} \alpha_{k} d\alpha_{k+1} \ldots d\alpha_n\\
&+&\sum_{k=1}^{n-1}
(-1)^{\sum_{j=0}^{k}(n-j-1)\partial \alpha_j
+ \sum_{j=k+1}^{n}(n-j)\partial \alpha_j + k}
\alpha_0d\alpha_1\ldots d\alpha_{k} \alpha_{k+1} d\alpha_{k+2}\ldots d\alpha_n\\
&+&
(-1)^{\nu + (n-1)\partial \alpha_n\sum_{j=0}^{n-1}(n-j-1)\partial \alpha_j + n}
\alpha_n\alpha_0d\alpha_1\ldots d\alpha_{n-1}
\end{eqnarray*}
The first term cancels with the first term in the first sum, the second term in the
first sum cancels with the first  term in the second sum, and so on. The only terms 
that do not cancel are the last term in the second sum and the last term.
The computation then follows.
$\;\;\sqcap \!\!\!\! \sqcup$ \end{proof}
\begin{corollary}
\label{corollary1}
 Let $\Omega=F_0\Omega\supset F_{-1}\Omega
\supset \ldots$ be a filtered differential superalgebra,
$\tau : \Omega \to \Bbb C$ a supertrace satisfying 
$\tau (F_{-n-1}\Omega)=0$ and $\rho:A=\Omega/F_{-1}\Omega \to \Omega$
be an arbitrary set theoretic lifting.  Define $\mu$ as in lemma \ref{lemma1}.
Then $\varphi(a_0,a_1,\ldots,a_n)=
(-1)^\mu \tau(\rho(a_0)d\rho(a_1)\ldots d\rho(a_n))$
is a Hochschild cocycle independent on the lifting $\rho$.
\end{corollary}
\begin{proof}
Define $\tilde{\varphi}:\Omega^{\otimes n+1}\to \Bbb C$, $\tilde{\varphi}=
\tau\circ \Phi$. Then the previous lemma shows that
$\tilde{\varphi}\circ b=0$. Moreover if any of $\alpha_0,\ldots,\alpha_n$
is in $F_{-1}\Omega$ then $\tilde{\varphi}(\alpha_0,\ldots,\alpha_n)=0$.
The lemma then follows.
$\;\;\sqcap \!\!\!\! \sqcup$ \end{proof} 
\begin{lemma}
\label{lemma2}
If in the above Corollary one also has $\tau(da)=0$ for any $a\in \Omega$
then $\varphi$ has the cyclic property (equation (\ref{eq3})).
\end{lemma}
\begin{proof} We have by definition
\begin{eqnarray*}
\lefteqn{\varphi(a_n,a_0,\ldots,a_{n-1})}\\
&=& (-1)^{\partial a_{n-2}+\partial a_{n-4}
+\ldots +n\partial a_n}\tau(a_nda_0da_1\ldots da_{n-1})\\
&=& (-1)^{\partial a_{n-2}+\partial a_{n-4}
+\ldots +\partial a_n(\partial a_0 +\ldots
\partial a_{n-1})}\tau(da_0da_1\ldots (da_{n-1})a_n)\\
&=& \!(-1)^{\partial a_{n-2}+\partial a_{n-4}
+\ldots +\partial a_n(\partial a_0 +\ldots
\partial a_{n-1})+
\partial a_0 +\partial a_1 + \ldots \partial a_{n-1} +n}
\tau(a_0da_1\ldots da_n)\\
&=& (-1)^{\nu+n}\varphi(a_0,a_1,\ldots,a_n)
\end{eqnarray*}
where $\nu$ is as in equation (\ref{eq3}).
$\;\;\sqcap \!\!\!\! \sqcup$ \end{proof}

We now review an important construction due to Connes.

Let  $(\Omega,\nabla,\tau)$ be a cycle with curvature
$\omega$.
Define $(\tilde{\Omega},d,\tilde{\tau})$ by
\begin{equation}
\label{eqOmega} 
\tilde{\Omega}=\Omega\oplus \Omega X\oplus X\Omega \oplus X\Omega X
\end{equation}
 where $X$ 
is a formal odd degree symbol and is not to be
considered alone (i.e. $aX$, $Xa$ do make sense, but $X$ does not). 
The filtration on $\tilde{\Omega}$ is defined by
$F_{-k}\tilde{\Omega}=F_{-k}\Omega\oplus F_{-k+1}\Omega X\oplus 
XF_{-k+1}\Omega \oplus XF_{-k+2}\Omega X$.
On $\tilde{\Omega}$ we consider the 
multiplication 
$(aX)b=a(Xb)=0\;\; \forall a,b \in \Omega$, $(aX)(Xb)=a\omega b$, 
and the differential
\begin{equation}
\label{eq1}
da=\nabla a +Xa+ (-1)^{\partial a}aX\;\;\;\forall a\in \Omega, \;\;\;\;dX=0 
\end{equation}
The trace is  $\tilde{\tau}(a_{00}+a_{01}X+Xa_{10}+Xa_{11}X)=
\tau(a_{00})-(-1)^{\partial a_{11}}\tau(\omega a_{11})$.

\begin{theorem}[Connes] 
\label{thm1}
(i)  Using the above notation 
$(\tilde{\Omega},d,\tilde{\tau})$ is an exact cycle. 
If $\tau(F_{-n-1}\Omega)=0$ then $\tilde{\tau}(F_{-n-1}\tilde{\Omega})=0$.
Moreover we have  $\Omega=e\tilde{\Omega}e$ and  $\nabla(a)=e(da)e$ if 
$a\in \Omega$ and $\Omega$ has a unit $e$.

(ii) Suppose $\Omega$ is a filtered  algebra endowed with a connection $\nabla$ and 
curvature $\omega$. 
Let $\rho:B\to\Omega$ be a degree preserving algebra morphism.
Suppose $B$ has a unit denoted $e$.
Then there exists a morphism $\varphi:e\Om(B)e\to \Omega$ such that
$\nabla(\varphi(b))=\varphi(e(db)e)$ and $\omega=\varphi(edede)$.
$\;\;\sqcap \!\!\!\! \sqcup$ 
\end{theorem}

Note that the sign in (\ref{eq1}) is different from 
the one in \cite{Co2}, page 329. 
This is necessary in order to have $d(aXb)=0$ and 
for $\tilde{\tau}$ to be a closed trace,
as seen from the computation of $\tilde{\tau}(d(aX))$. 

It makes sense to adjoin $X$ to the algebra but then 
$X^2\not = \omega$ since $d(X^2)=0=\nabla \omega \not = d\omega$. This 
would not contradict $(aX)X=a\omega$ since the unit $e$ of $\Omega$, 
if it exists,
is never a unit of the bigger algebra: $eXa\not =Xa$ if $a\not =0$.

Part (ii) is the analog in the  connection case of the corresponding
result for differential graded algebras and immediately follows from that
case using (i). Indeed let $\tilde{\Omega}=\Omega\oplus \Omega X\oplus X\Omega
\oplus X\Omega X$
be the differential algebra associated to $(\Omega,\nabla)$ as in (i).
Then by universality  we get  a morphism  $\Om(B) \to
\tilde{\Omega}$
whose restriction  $e\Omega^*(B)e \to \Omega =e\tilde{\Omega}e$ is the
morphism $\varphi$ we are looking for. The last statement follows
from $\varphi(edede)=e(\nabla e +Xe +eX)(\nabla e +Xe +eX)=\omega$ since
$\nabla e=0$.
This is in agreement with equation (\ref{eq1}).
\begin{theorem}[Connes] 
\label{thm2}
Let  $(\Omega,d,\tau)$ be an exact cycle,  and define
$A=\Omega/F_{-1}\Omega$. 
Suppose that $\tau(F_{-n-1}\Omega)=0$ and choose an arbitrary set
theoretic lifting $\rho:A\to \Omega$ for $\Omega\to A$. Then 
$$\varphi(a_0,a_1,\ldots,a_n)=
(-1)^\mu \tau(\rho(a_0)d\rho(a_1)\ldots d\rho(a_n))$$
is a cyclic cocycle on $A$ independent of the choice of
$\rho$. Here $\mu=\partial a_{n-1}+\partial a_{n-3}+
\partial a_{n-5}+\ldots$.
$\;\;\sqcap \!\!\!\! \sqcup$ 
\end{theorem}
The above theorem gives a (non-commutative) geometric way of
defining cyclic cocyles and was considered by Connes for trivially
graded algebras in \cite{Co2}. Together with Theorem 
\ref{thm1} it provides us with a canonical method of
construction cyclic cocycles. The proof is contained in the lemmata above.
The extra sign $\mu$ is necessary in order
to get the Hochschild cocyle property in the $\Bbb Z_2$-graded case. 
It is interesting to note that the sign was originally obtained
from the formula of
the Fedosov product using the bivariant Chern-Connes character
\cite{Ni4}.
A different formula was considered
in \cite{Kastler}. 

Assuming the trace 
$\tau$ to be even, i.e. $\tau(a)=0$ if $a$ is odd, we also
obtain using \cite{NiD} an odd  cyclic
cocycle $\psi^\tau$ on the crossed product algebra
$A\rtimes \Bbb Z_2$ defined by 
\begin{eqnarray}
\label{eq6-}
\lefteqn{2\psi^\tau(a_0v^{i_0},a_1v^{i_1},\ldots, a_nv^{i_n})=} \nonumber \\
\!\!\!&=&\!\!\!
\left \{
\begin{array}{ll}
(-1)^{\nu} \varphi^\tau(a_0,\ldots,a_n)=(-1)^{\mu+\nu} 
\tau(a_0da_1\ldots da_n)
\\
\mbox{ if }\partial a_0+\ldots+\partial a_n
\mbox{ is even and } 
i_0+\ldots+i_n \mbox{ is odd } \\
0 \mbox{ otherwise. }
\end{array}
\right.
\end{eqnarray}
Here 
$\nu=0$ if $n=0$ and 
$\nu=\sum_{k<n}i_k(\partial a_{k+1}+\ldots+\partial a_{n})$
if $n>0$ and $\mu$ is as in the above lemma. Moreover $v^2=1$ is the
invertible element implementing the action of 
$\Bbb Z_2$: $vav^{-1}=\alpha(a)$.
In the notation of \cite{NiD} we have 
$\psi^\tau=p_{A}^*{\varphi^\tau}$.

\begin{definition}
\label{def2.7.}
Let  be a cycle with the property that
$\tau(F_{-n-1}\Omega)=0$. The {\em cyclic cocycle associated to 
$(\Omega,\nabla,\tau)$} (and $n$) is by definition the cyclic cocycle
$\varphi^\tau:A^{\otimes n+1}\to \Bbb C$ defined using Theorems
\ref{thm1} (i) and \ref{thm2}.
\end{definition}

\paragraph{ } One of the reasons for working in the setting of 
filtered algebras rather than that of differential graded algebras,
besides the connections with the Fedosov product, is that
it is closer to the original idea of Quillen \cite{Quillen0}
of using higher homogeneous components in the superconnection, in addition 
to the usual degree one component.
This is crucial for the following theorem. 

Denote
by $HC^n(A)$ the cyclic cohomology groups of a superalgebra $A$ as
defined in \cite{Co2} with the necessary changes required by 
the grading (see \cite{Kassel3}). Also let  $S:HC^{n}(A)\to HC^{n+2}(A)$ 
be Connes' periodicity operator \cite{Co2}.

The proof of the following theorm will require the notion of cobordism
\cite{Co2}. 

\begin{definition}[Connes] 
\label{def-co}
Two cycles $(\Omega_0,\nabla_0,\tau_0)$ and 
$(\Omega_1,\nabla_1,\tau_1)$ will be  called {\em cobordant} if there exists a 
chain $(\Omega',\nabla',\sigma)$ and a morphism 
$r:\Omega'\to \Omega_0\oplus\Omega_1$
which is compatible with connections and curvatures
\begin{equation}
\label{eq-comp}
r(\nabla'a)=(\nabla_0\oplus\nabla_1)r(a), \;\;r(\omega')=\omega_0\oplus\omega_1
\end{equation} 
(we implicitely assume that $r$ can be extended to include
$\omega'$ in its range) and which satisfies Stoke's Theorem
\begin{equation}
\label{Stokes}
\sigma(\nabla'a)=(-\tau_0\oplus\tau_1)(r(a))
\end{equation}
\end{definition}

Connes' main result on cobordant cycles is that two such {\em closed}
cycles give rise to cyclic cocycles $\varphi_0,\varphi_1$ satisfying 
$S\varphi_0=S\varphi_1$ \cite{Co2} Theorem 32.
(The actual theorem states that $\varphi_1-\varphi_0=B\psi$ for some 
Hochschild cocyle $\psi$, which is equivalent to the statement we have
given above using Connes' exact sequence, loc. cit.) In the following
we shall sketch the easy extension of this result to arbitrary cycles
(i.e. not necessarily satisfying $\nabla^2=0$).
\begin{lemma}
(i) The functor $\;\tilde{ }$, 
$\Omega \to \tilde{\Omega}$, from cycles to 
closed cycles defined in theorem \ref{thm1} preserves direct
sums and is linear in $\tau$.

(ii) Moreover the functor $\;\tilde{ }$ preserves the cobordism between cycles.
\end{lemma}
\begin{proof} The first part is immediate from the definitions.

In order to prove (ii) we shall use the notation in Definition
\ref{def-co} and Theorem \ref{thm1} (i). The compatibility equation
(\ref{eq-comp}) is obvious. The Stokes theorem (\ref{Stokes}) is 
immediate for $a\in \tilde{\Omega'}$ of the form $\alpha, \alpha X$ or 
$X\alpha$ with $\alpha \in \Omega'$. Suppose now that $a=X\alpha X$ 
with $\alpha \in \Omega'$, then we have
\begin{eqnarray*}
\lefteqn{\tilde{\sigma}(d(X\alpha X))=}\\
&=&-(-1)^{\partial \alpha}\sigma(\omega' {\nabla}'\alpha)\\
&=&-(-1)^{\partial \alpha}\sigma({\nabla}'(\omega'\alpha))\\
&=&-(-1)^{\partial \alpha}(-\tau_0)\oplus \tau_1(r(\omega' \alpha))\\
&=&-(-1)^{\partial \alpha}(-\tau_0)\oplus \tau_1
((\omega_0\oplus\omega_1)r(\alpha))\\
&=&(-\tilde{\tau}_0)\oplus \tilde{\tau}_1(Xr(\alpha)X)
\end{eqnarray*}
This proves the lemma.
$\;\;\sqcap \!\!\!\! \sqcup$ \end{proof}
\begin{corollary}[Connes]
Suppose the cycles $(\Omega_0,\nabla_0,\tau_0)$ and
$(\Omega_1,\nabla_1,\tau_1)$ are cobordant and satisfy
$\tau_i(F_{-n-1}\Omega_i)=0$. Let $\varphi_0$ and $\varphi_1$ the
cyclic cocyles associated to these cycles. Then $S\varphi_1=S\varphi_0$.\sqep 
\end{corollary}
\begin{theorem}
Let $(\Omega,\nabla,\tau)$ be a cycle.
Choose an arbitrary odd element $\eta \in F_{-1}\Omega$. Define
$\nabla_1=\nabla+ad_\eta$ and
$\varphi_0$, respectively $\varphi_1$, to be the cyclic cocycles associated
to $(\Omega,\nabla_0=\nabla,\tau)$, respectively to $(\Omega,\nabla_1,\tau)$.
Then $S\varphi_0=S\varphi_1$.
\end{theorem}
\begin{proof} The proof will follow from 
the above corollary
once we prove that the 
two  cycles  in the theorem are cobordant.
The required cobordism is given by $(\Omega[0,1]+\Omega[0,1]dt,\nabla'=
\nabla+dt\frac{\partial}{\partial t}+ad_{t\eta},\sigma)$
where $\sigma(a+(dt)b)=\int_0^1dt\tau(b(t))$. The morphism
$r$ is the evaluation at the end points. 
(Note that the physicist's notation for integrals is the correct one
in the case of noncommuting variables!) 
The fact that $\sigma$ is a trace and the relation
$\sigma(\nabla_1(a))=\tau(a(1)-\tau(a(0))$ are easy computations. 
$\;\;\sqcap \!\!\!\! \sqcup$ \end{proof}

\section{Topological preliminaries}
In this section we shall discuss some aspects related to the 
topologies on the algebras that we are going to work with.
We also prove the single generation for the homology of a certain 
complex similar to the complexes used in \cite{NiD}.

We shall frequently use results and 
definitions from  \cite{NiD} and we 
begin by recalling some of them.

Let ${\cal D}_{\epsilon}$ ($\epsilon>0$) be the subalgebra 
of $C^{\infty}(\Bbb R)$
of those functions satisfying the following exponential decay condition
at $\infty$
$$p_{k,\epsilon}(f)^2=\int_{\Bbb R} 
e^{\epsilon x^2}|f^{(k)}(x)|^2dx<\infty \;\;(\forall)
k\geq 0$$
endowed with the  Frechet topology generated by the above family
of seminorms. Then let ${\cal D}=\displaystyle{\lim_{\rightarrow}} 
\;{\cal D}_{1/n}$ with the
corresponding inductive limit topology.
This definition of ${\cal D}$ coincides with the one in \cite{NiD} although
the ${\cal D}_\epsilon$ are different. The new definition has the advantage that
it easely gives the following result.
\begin{proposition}
\label{prop-nucl} 
The spaces ${\cal D}_\epsilon$ ($\epsilon >0$) are nuclear,
and hence ${\cal D}$ is also nuclear.
\end{proposition}
\begin{proof} The last part is an immediate consequence of the properties
of the inductive limit \cite{Treves}.

Fix $\epsilon >0$. 
It is enough to show that for each of the defining seminorms 
$p=p_{k,\epsilon}$ there exists an other seminorm $q$ such that the 
inclusion of Hilbert spaces $({\cal D}_\epsilon)_q\to  
({\cal D}_\epsilon)_p$ is nuclear \cite{Treves} (recall that ``nuclear'' means
``trace class'' when only Hilbert spaces are involved). 
We can also restrict to the subspace of the functions vanishing on
$(-\infty,0]$.
If one takes $q=p_{k+2,\epsilon}$ then the above spaces are naturally
isometrically isomorphic to $L^2( [0,\infty))$ with the standard Lebesque
measure and the inclusion becomes the integral operator with
kernel $K(x,y)=0$ if $x>y\geq 0$ and $K(x,y)=(-x+y)
\exp(\epsilon(x^2-y^2)/2)$ if $0\leq x \leq y$. Since this kernel is square
integrable the corresponding operator is Schatten-von Newmann 
(i.e. $Tr(A^*A)<\infty$). The product of two Schatten-von Newmann 
operators is trace class and hence we get the result if we let 
$q=p_{k+4,\epsilon}$.$\;\;\sqcap \!\!\!\! \sqcup$ 
\end{proof} 

Denote by $\Sigma^{(n)}$ the Fr\'{e}chet space of symbols of
order $n$ on $\Bbb R$:
\begin{equation}
\!\!\Sigma^{(n)}\!=\!\{f:\Bbb R\to \Bbb C: 
\,\forall k\geq 0 \, \exists C_k>0\;\;
|f^{(k)}(x)|<C_k(1+|x|)^{n-k} \;\; \forall x\in \Bbb R\! \}
\end{equation}
The intersection $\Sigma^{(-\infty)}=\cap \Sigma^{(n)}$ 
is the Schwartz space of 
rapidly decreasing functions on $\Bbb R$, also denoted 
by ${\cal S}$, it has a natural 
Fr\'{e}chet topology. The union $\Sigma^{(\infty)}=\cup \Sigma^{(n)}$
of all symbols will have the inductive limit topology. 
We endow ${\cal S}$ and $\Sigma^{(\infty)}$ with the grading
$\alpha(f)(x)=f(-x)$.
 
Let $\widehat{\otimes}$ denote the (completed) projective tensor product 
\cite{Grothendieck} and consider  the spaces 
$$
{\cal F}_k(m_0,\ldots,m_n)
=\Sigma^{(m_0)}\widehat{\otimes} 
\ldots \widehat{\otimes} 
\Sigma^{(m_{k-1})} \widehat{\otimes}  {\cal D}_{1/m_k} \widehat{\otimes}  
\Sigma^{(m_{k+1})} \widehat{\otimes}  
\ldots \widehat{\otimes} \Sigma^{(m_n)}
$$
where $m_0,\ldots m_n \in \Bbb N$.
Also define 
\begin{equation}
\label{missing-eq}
{\cal F}_k^{(n)}=\lim_{\to }{\cal F}(m_0,m_{1},
\ldots,m_n)
\end{equation}
where in the inductive limit all $m_j$ go to $\infty$. Obviously
${\cal F}_k^{(n)}$ contains $\Sigma^{(\infty)\otimes k-1}
\otimes {\cal D} \otimes \Sigma^{(\infty)\otimes n-k}$.
We want to consider the union of all the spaces ${\cal F}_k^{(n)}$.
\begin{proposition}
\label{prop2.2}
The natural map ${\cal F}_k^{(n)}\to C^\infty(\Bbb R^{n+1})$ is injective.
\end{proposition}
The idea of the proof is to reduce everything to the compact case using the
following lemma. Denote by $l^1$ the Banach space of absolutely summable
sequences of complex numbers.
\begin{lemma}
\label{lemma2.4}
Fix $a>1$ and let $\chi:\Bbb R \to [0,1]$ be an even smooth function with support in $[-a,a]$
satisfying $\chi=1$ on $[-1,1]$. Denote by $\chi_0=\chi$ and 
$\chi_n(x)=\chi(x/a^{2n})-\chi(x/a^{2(n-1)})$ for $n\geq 1$.
Then for any $m$ the map $T:\Sigma^{(m)}\to l^1 \widehat{\otimes} \Sigma^{(m+1)}$,
$$T(f)=(\chi_0f,\chi_1f,\ldots,\chi_nf,\ldots)$$
is continuous and in particular $f=\sum \chi_n f$ absolutely in $\Sigma^{(m+1)}$ for
any $f\in \Sigma^{(m)}$. The same result holds true for 
${\cal D}_{1/m}\to l^1\widehat{\otimes}{\cal D}_{1/(m+1)}$.
\end{lemma} 
\begin{proof}
Indeed since the multiplication 
$\Sigma^{(p)}\otimes \Sigma^{(q)} \to \Sigma^{(p+q)}$
is continuous it is enough to prove that $1=\sum \chi_n$ absolutely in
$\Sigma^{(1)}$. To show this observe first that due to the definition of
$\chi$ the support of $\chi_n$ is contained in $\{|x|\leq a^{2n+1}\}$ for $n\geq 1$.
Hence using $\chi_n(x)=\chi_1(x/a^{2(n-1)})$ we get
$$\sup_x
(1+|x|)^{p-1}|\chi_n^{(p)}(x)|\leq C \sup_x |x|^{p-1}|\chi_n^{(p)}(x)|
\leq C a^{-2(n-1)} \sup_y |y|^{p-1}|\chi_1^{(p)}(y)|$$
where $C$ is a  positive constant independent on 
$n$ and $y=x/a^{2(n-1)}$.
The estimates for $|1-\sum_0^N\chi_n|$ are obtained 
using the same technique.
$\;\;\sqcap \!\!\!\! \sqcup$ \end{proof} 
We now go back to the proof of the Proposition \ref{prop2.2}
\paragraph{Proof of the proposition.}
Suppose $f\in {\cal F}_k(m_0,\ldots,m_n)$ maps to $0$ in $C^\infty(\Bbb R^{n+1})$.
Let, for any  multiindex $\alpha=(k_0,\ldots,k_n)$, $\chi_\alpha=
\chi_{k_0}\chi_{k_1}\ldots\chi_{k_n}$. Also 
let $|\alpha|=\max\{k_0,\ldots,k_n\}$ for
$\alpha$ as above. Our assumption on $f$ tells us that each of $f\chi_\alpha$ maps
to $0$ in $C^\infty(\Bbb R^{n+1})$. Since both $C_c^\infty(\Bbb R)$ and
$C^\infty(\Bbb R)$ are nuclear we can replace the projective tensor product by
the injective tensor product and the latter preserves the injectiveness. 
This shows that $f\chi_\alpha=0$ in ${\cal F}_k(m_0,\ldots,m_n)$.
Since $f=\sum_\alpha f\chi_\alpha$ in ${\cal F}_k(m_0+1,\ldots,m_n+1)$
thanks to the above lemma we obtain that $f$ maps to $0$ in 
${\cal F}_k(m_0+1,\ldots,m_n+1)$. This proves the proposition.
$\;\;\sqcap \!\!\!\! \sqcup$ 

Using the proposition we have just proved we may define 
\begin{equation}
\label{eq-def}
{\cal F}^{(n)}=\sum_{k} {\cal F}_k^{(n)}
\end{equation}
as a subspace of $C^\infty(\Bbb R^{n+1})$.

From the definition it follows that 
$$
{\cal D}^{\widehat{\otimes} n+1}
\subset {\cal F}^{(n)}\subset 
C^\infty(\Bbb R)^{\widehat{\otimes} n+1}=C^\infty
(\Bbb R^{n+1}).
$$ Since $\Sigma^{(p)}\otimes {\cal D}_{1/p}\to
{\cal D}_{1/(p+1)}$ is continuous the Hochschild differential $b$ 
maps ${\cal F}^{(n)}$ to ${\cal F}^{(n-1)}$ and 
the above inclusions are chain maps. 

We now begin the computation 
of the homology of $({\cal F}^{(n)},b)$ which will prove useful
in the next section.

Recall that this Hochschild boundary is defined using the grading,
see equation (\ref{eq2}). In addition to this boundary which we will
also denote in the following by $b_{grad}$ in order to avoid any
confusion, we will also use $b_{ev}$ and $b_{odd}$.
The second formula is that of the  Hochschild boundary defined
ignoring of the grading (i.e. consider all $\partial a=0$). The formula
for $b_{ev}$ is given by 
\begin{eqnarray*}
b_{ev}(\Tt)&\!=\!&a_0a_1\otimes\ldots\otimes a_n-a_0\otimes 
a_1a_2\otimes\ldots\otimes a_n+\ldots\\
&\!-\!&(-1)^{n}a_0\otimes\ldots\otimes a_{n-1}a_n +(-1)^n \alpha(a_n)a_0\otimes
\ldots\otimes a_{n-1}
\end{eqnarray*}
In the above formula $\alpha$ is the grading morphism.

As seen from the definition these differentials differ only in the
form  of the last term. Moreover they all commute with the
action of the grading.
Actually the main reason for considering all 
these differentials is that $b_{grad}=b_{ev}$ on the {\em even} parts of
the complexes, and $b_{grad}=b_{odd}$ on the {\em odd} parts. 
See also \cite{NiD}, Proposition 1. These identification will 
simplify the computation of the $b_{grad}$ homology as in \cite{NiD}.

We will discuss in detail the $b_{ev}$ homology which is actually 
the result needed in the next section. The proof in the odd case is 
similar. We obtain that
the extra unbounded operators that we introduce do not change the
homology which is again singly generated \cite{NiD}.

%

Denote $\chi_n^{[m]}=H(nx)(\chi^{[m]}(2^{-m-n/m}x)-\chi^{[m]}(2^{-m-(n-1)/m}x))$
where $\chi^{[m]}:\Bbb R\to [0,1]$ is a smooth function with support in
$[-2^{1/(2m)},2^{1/(2m)}]$, $\chi^{[m]}=1$ on $[-1,1]$, and $H$ is the 
Heaviside function, $H(x)=1$ if $x\geq 0$, $H=0$ otherwise.

Let $\chi_{\alpha}^{[m]}=\chi_{\alpha_0}^{[m]}
\otimes 
\ldots \otimes \chi_{\alpha_n}^{[m]}$ for any multiindex $\alpha=(\alpha_0,
\ldots,\alpha_n)\in \Bbb Z^{n+1}$.
Lemma \ref{lemma2.4} shows that $f=\sum_\alpha \chi_\alpha^{[m]} f$ for any 
$f\in {\cal F}^{(n)}$. Let us introduce the submodules 

$$
F_{i,m}{\cal F}^{(n)}=\{f\in {\cal F}^{(n)}, 
f=\sum_{\alpha \in J_{n,i}} \chi_\alpha^{[m]} f_0\}
$$
where $J_{n,i}=\{\alpha \in \Bbb Z^{n+1}, \mbox{ at least one of } 
|\alpha_n+\alpha_0|, |\alpha_0-\alpha_1|,
\ldots, |\alpha_{i-1}-\alpha_i| \mbox{ is } \geq 2\}$.
The first $+$ is not a mistake and is related to the 
twist in the 
definition of the $b_{ev}$. Also observe that the complement of
$J_{n,n}$ is a finite set.

The following lemma is the analog of Lemma 3.1 in \cite{Bryl-Nistor}
with the changes imposed by the topology. The proof will be similar.
\begin{lemma} 
The spaces $F_{i,m}{\cal F}^{(n)}$ 
are invariant under
$b_{ev}$, and for any
$i$ the complexes
$\cup_m F_{i,m}{\cal F}^{(n)}$ have vanishing homologies.
\end{lemma}
\begin{proof} The invariance follows from $\chi_p^{[m]}\chi_q^{[m]}=0$ if
$|p-q|\geq 2$. Let $\sigma:F_{i,m-2}{\cal F}^{(n)} \to F_{i,m}{\cal F}^{(n+1)}$
be given by 
$$\sigma(f)(x_0,\ldots,x_{n+1})
=\!\sum_{|p-q|\geq 2}\chi_p^{[m]}(x_i)\chi_q^{[m]}(x_{i+1})
f(x_0,\ldots,x_{i-1},x_{i+1},\ldots,x_{n+1})
$$
The formula $(b_{ev}\sigma +\sigma b_{ev} 
+(-1)^i)f\in F_{i-1,m-2}{\cal F}^{(n)}$
proves the result.
$\;\;\sqcap \!\!\!\! \sqcup$  
\end{proof}
\begin{proposition}
The $b_{ev}$ homology is singly generated concentrated in dimension
$0$ and the generator is even.
\end{proposition}
\begin{proof}
Define $F_\infty{\cal F}^{(n)}=\cup_{i,m} F_{i,m}{\cal F}^{(n)}$ and 
observe that it consists of those functions that vanish in  a neighborhood of the
$0$. This shows that 
the complex $({\cal F}^{(n)}/F_\infty{\cal F}^{(n)},b_{ev})$ 
has the same homology as $({\cal F}^{(n)},b_{ev})$. 
The first complex is ``independent of the topology'' since 
it is concentrated
at $0$ (this is the analog of the localization used in \cite{NiD})
and the same computation as in the above lemma shows that it 
also computes the homology of 
$(C^\infty([-1,1])^{\widehat{\otimes} n+1},b_{ev})$ which is
concentrated in dimension $0$ and is even as shown in \cite{NiD}.\sqep
\end{proof}
The complex ${\cal F}^{(n)}$ inherits also  actions of the cyclic
group \cite{Co3}: 
$$
t(\Tt)=(-1)^{n+\nu}a_n\otimes a_0 \otimes \ldots \otimes
a_{n-1}
$$ 
if $\nu=\partial a_n(\partial a_0+
\ldots \partial a_{n-1})
$ 
and one can define Connes' cyclic complex  
$({\cal F}^{(n)}/(1-t){\cal F}^{(n)},b_{grad})$ \cite{Co2}.
\begin{theorem}
\label{thm3.6}
The homology of the complex
$({\cal F}^{(n)}/(1-t){\cal F}^{(n)},b_{grad})$
has dimension one in each positive dimension and the generator
$c_n$ in dimension $n$ satisfies 
$\alpha(c_n)=(-1)^n c_n$. In dimension $0$
the even part of the homology is again singly generated by 
$c_0$. The generators can be chosen so that $Sc_{n+2}=c_n$.
\end{theorem}
\begin{proof} 
The proof now continues exactly as in \cite{NiD}.
One  sets the analog of the $(b,B)$-bicomplex
of Loday and Quillen from which one gets the Connes' long exact
sequence relating the cyclic homology (as defined 
by factoring $(1-t)$ as
above) and Hochschild homology. The result then  follows.\sqep
\end{proof}

\section{The Chern character of the index}
In this section we will be concerned with finding explicit formule
for the characteristic
numbers of the index. We have given such formulae in \cite{NiD}. 
In this section we show how super-connections provide an other
sequence of such formulae. This gives an other proof
of one of the basic formulae in \cite{Bismut2} Theorem 2.6 page 109,
see also \cite{BGV} Chapter 9. 

We endow ${\cal D}$ with the grading $\alpha(f)(x)=f(-x)$. 
Recall the followin definition from \cite{NiD}:
\begin{definition} A {\em generalized 
$\theta$-summable Fredholm operator} in $A$ 
is a  degree preserving continuous 
morphism $D:{\cal D} \to A$, where $A$ is a topological 
super-algebra. 
\end{definition}

A simple but crucial observation is that an operator $D$ as
in equation (\ref{eqD}) determines uniquely a generalized 
$\theta$-summable Fredholm operator (denoted by the same letter) 
in the algebra of bounded 
linear operators on a Hilbert space ${\cal H}$ via functional
calculus. Moreover, most importantly,
the graph projection (defined by equation (\ref{eqgraph}))
can be written down only in terms 
of $D$ and the grading automorphism $v$:
\begin{equation}
\label{eqind}
p=(1+v)/2+{D}(-e^{-x^2})v+D(\tau(x)x)
\end{equation}
where $\tau$ is a smooth  even  function satisfying
$\tau(x)^2x^2=e^{-x^2}(1-e^{-x^2})$ (as in introduction).

Motivated by this we define the index in general 
\begin{definition}
\label{def4.2}
The index of a $\theta$-summable Fredholm operator $D:{\cal D}\to A$ 
is by definition 
\begin{equation}
\label{eqIND}
\mbox{Ind}(D)=[e_0]-[p]
\in K_0(A\rtimes_{\alpha}\Bbb Z_2)
\end{equation}
where $e_0=(1+v)/2$ and $[e_0],[p]$ denote the class of the 
idempotents $e_0$ and $p$ in the corresponding $K$-theory group.
\end{definition}

Note that $\mbox{Ind}(D)$ is an odd class in 
$K_0(A\rtimes \Bbb Z_2)$ with respect to the `dual' action
of $\Bbb Z_2$ on $A\rtimes \Bbb Z_2$ given by 
$\hat{\alpha}(a+bv)=a-bv$.

Consider the space
\begin{equation}
\label{eqcalD}
\overline{\Omega}=\prod_{n=0}^\infty \overline{\Omega}^n
\, , \,
\overline{\Omega}^0={\cal F}^{(0)}
\, , \, 
\overline{\Omega}^n=\Bbb C1\!\otimes
{\cal F}^{(n-1)} \oplus {\cal F}^{(n)},n>0
\end{equation}  
with differential given by 
$d(m)=1\otimes m$, $d(1\otimes m)=0$. Obviously $d^2=0$. 
We see from the definition and proposition \ref{prop2.2}
that we have the following inclusions 
$$\Omega^*({\cal D}) \subset 
\overline{\Omega} \subset \prod_n\Omega^n(C^\infty(\Bbb R))$$
Using Lemma \ref{lemma1} and the fact that ${\cal F^{(n)}}$ is
invariant under the $b$ we see that the multiplication
of $\Om(C^\infty(\Bbb R))$ restricts to a multiplication of 
$\overline{\Omega}$. Similarily the left and right representations
of $\Sigma^{(\infty)}=\cup \Sigma^{(n)}$  on 
$\Om(C^\infty(\Bbb R))$ restrict to a
representations of $\Sigma^{(\infty)}$ on $\overline{\Omega}$
which are  always nonunital 
since $\Sigma^{(\infty)}(\Bbb C1\!\otimes {\cal F}^{(n-1)})
\subset {\cal F}^{(n)}$. These representations are 
determined by the fact that the algebra
$\overline{\Omega}$ has as a ``dense'' subset the linear
span of $a_0da_1\ldots da_n$ and
$da_0da_1\ldots da_n$ where $a_0,\ldots,a_n \in \Sigma^{(\infty)}$
and at least one of them is in ${\cal D}$. 

We are going to use the concept of a cycle
introduced in definition 2.1. (iii).
\begin{definition}
\label{defFS} A {\em finite summable cycle} is a
cycle $(\Omega,\nabla,\tau)$ toghether with a morphism 
$\rho:{\cal D}\to \Omega$ satisfying:

(i) the natural morphism 
$\Om({\cal D})\to \tilde{\Omega}$ defined by $D$
extends to a  morphism
$\overline{\rho}:\overline{\Omega}\to \tilde{\Omega}$,
and 

(ii) $\tilde{\tau}$ vanishes on $\overline{\rho}([\Sigma^{(\infty)},
\overline{\Omega}])$.
\end{definition}
In the above definition  $\overline{\Omega}\supset\Om({\cal D})$ is as 
defined above, and 
$\tilde{\Omega}=\Omega\oplus 
\Omega X\oplus X\Omega \oplus X\Omega X$ and $\tilde{\tau}$ (the 
extension of $\tau$) are 
as in Theorem  \ref{thm1}. 

\paragraph{Remark.} We also have 
$\tilde{\tau}\circ\rho([d\Sigma^{(\infty)},\overline{\Omega}])=0$. 
This is obtained using the formula $[df,g]=d[f,g]-(-1)^{\partial f}[f,dg]$.

Let ${\frak B}$ be  a filtered topological superalgebra containing $\Omega$ as an ideal
so that 
$F_{-n}\Omega =\Omega \cap F_{-n}{\frak B}$, and
the inclusion  $\Omega\to {\frak B}$ and the multiplications 
$\Omega\times {\frak B}\to \Omega$ and ${\frak B}\times \Omega\to \Omega$ 
are continuous.
Also assume that the connection $\nabla$ on
$\Omega$ extends to a continuous connection, also denoted $\nabla$, 
on ${\frak B}$, and that $F_{-N}{\frak B}=0$ for some large
$N$. The following Proposition gives a wide class of  examples of
finite summable cycles.
\begin{proposition}
\label{prop4.4} 
Let $(\Omega,\nabla,\tau)$ be a cycle, 
$D:{\cal D}\to \Om$ be a 
continuous morphism, and ${\frak B}$ be as above. Also assume that 
there exists a continuous extension of $D$ to a morphism
$\overline{\rho}:\Sigma^{(\infty)} \to {\frak B}$ and that the trace $\tau$ satisfies 
$\tau([{\frak B},\Omega])=0$.
Then $(\Omega,\nabla,\tau)$ is a finite summable cycle.
\end{proposition}
\begin{proof} Consider
$$
D_{\#}:\Sigma^{(j_0)}\otimes 
\ldots \otimes \Sigma^{(j_{k-1})} \otimes {\cal D} 
\otimes \Sigma^{(j_{k+1})}\otimes 
\ldots \otimes \Sigma^{(j_n)}  \rightarrow  F_{-n}\tilde{\Omega}=\Omega\cap 
F_{-n}{\frak {\frak B}}
$$
$$
D_{\#}(\Tt)=\overline{\rho}(a_0)d\overline{\rho}(a_1)
\ldots d\overline{\rho}(a_n)
$$ 
Our continuity assumptions on $D$ imply  
that the above map extens to a continuous
map ${\cal F}^{(n)}\to F_{-n}\tilde{\Omega}$, 
where ${\cal F}^{(n)}$ is as defined in the previous 
section, equation (\ref{eq-def}).\sqep
\end{proof} 

The main reason for introducing the algebra $\overline{\Omega}$ and the
notion of $\theta$-summable cycles is to make sense of the exponentials
of the curvature of the super-connection
$\,D+\nabla$. This will then allow us 
to extend Quillen's method \cite{Quillen0}.

Denote by $\Delta_n=\{(t_0,\ldots,t_n)
\in \Bbb R^{n+1}: t_i\geq 0, t_0+\ldots+t_n=1\}$ the unit simplex in 
$\Bbb R^{n+1}$ and by $t\Delta_n$ its dilation by $t$. Also 
denote by $e$ the unit of $\Sigma^{(\infty)}$, such that
$e\overline{\Omega}e= {\cal F}^{(0)}\oplus \prod_{n=0}^\infty {\cal F}^{(n)}$
and define $\nabla(\xi)=e(d\xi)e$. We  introduce the
short hand notation
$\exp(-tD^2)$ for $D(f)\in \Omega$ if
$f(x)=\exp(-tx^2)$.

The following Lemmata are devoted to make sense in our setting of the  
well known  formal perturbative expansion:
\begin{equation}
\label{eq6}
\exp(-t(D+s\nabla)^2)=\sum_{n=0}^\infty\tilde{\eta}_n(s;t)
\end{equation}
where 
$$
\tilde{\eta}_n(s;t)= (-1)^n\int_{t\Delta_n} \exp(-t_0 D^2)
a_0 \exp(-t_1 D^2) a_0 \ldots
a_0 \exp(-t_n D^2)dV
$$ 
$a_0=s\nabla(\;D)+s^2\nabla^2$, and $dV$ is the volume element 
of $t\Delta_n$ normalized such that $Vol(t\Delta_n)=t^n/n!$.

We shall make sense of the above formal expression 
for finite summable cycles by defining 
$\eta_n(s;t) \in \overline{\Omega}$ which satisfies
$\overline{D}(\eta_n(s;t))=\tilde{\eta}_n(s;t)$, where $\overline{D}:
\overline{\Omega}\to \tilde{\Omega}$ 
is the extension of $D:\Om({\cal D})\to \tilde{\Omega}$.

Consider the functions 
$g_n:(0,\infty)\times [0,\infty)^{n+1} \to \Omega^*(C^\infty(\Bbb R))$, 
$g_0(s,t_0)=e^{-t_0x^2}$ and 
$g_n(s,t_0,\ldots,t_n)=e^{-t_0x^2}
a e^{-t_1x^2} a \ldots
a e^{-t_nx^2}$ where $x=id$ is the identity function of $\Bbb R$,
$a=se(dx)e+s^2edede=s(e\otimes x +x \otimes e) 
+s^2 (e\otimes e \otimes e) \in e\overline{\Omega}e$. Remember that $x$ is odd
and $d$ is a graded derivation
Also it is important to stress that 
the multiplication is the one induced by $\Omega^*(C^\infty(\Bbb R))$.
As a remark on notation we make the convention 
that  $e^{-tx^2}=e_0$ for $t=0$, where $e_0$ is the 
identity of $C^\infty(\Bbb R)$, different from  the identity of 
$\Om(C^\infty(\Bbb R))$. This is in agreement with  continuity requirements
(see Lemma \ref{l-cont} bellow). We shall write $e$ instead of $e_0$
bellow.
	
It follows from the definition that the functions  $g_n$
defined above are continuous and satisfy  
\begin{eqnarray}
\lefteqn{g_{n+1}(s;t_0,\ldots,t_{n+1})=}\nonumber\\
&=&
\label{eq9}
(sg_n(s;t_0,\ldots,t_n)dx+s^2g_n(s;t_0,\ldots,t_{n}) dede)
e^{-t_{n+1}x^2}\\
&=&
\label{eq10}
e^{-t_{0}x^2}
(sdxg_n(s;t_1,\ldots,t_{n+1})+s^2dedeg_n(s;t_1,\ldots,t_{n+1})) 
\end{eqnarray} 
\begin{lemma}
\label{l-cont} 
Define for any $t\geq 0$ the function $f_t(x)=e^{-tx^2}$.
Then the function 
$[0,\infty)\ni t\to f_t\in \Sigma^{(1)}$  is continuous.
\end{lemma}
\begin{proof}
This function is clearly continuous for $t>0$ so we only need
to check the continuity at $0$.

The continuity at $0$ is equivalent to 
\begin{equation}
\label{eq-cond1}
\lim_{t\to 0}\sup_{x\in \Bbb R}|1-e^{-tx^2}|(1+|x|)^{-1}=0
\end{equation}
and 
\begin{equation}
\label{eq-cond2}
\lim_{t\to 0}\sup_{x\in \Bbb R}|(e^{-tx^2})^{(n)}|(1+|x|)^{n-1}=0\;,\;\; n\geq 1
\end{equation}
The first condition (\ref{eq-cond1}) is easy so we concentrate on the
second one. We have
$$
\sup_{|x|\leq t^{-1/2}}|(e^{-tx^2})^{(n)}|(1+|x|)^{n-1}
\leq \sup_{|x|\leq 1}|(e^{-x^2})^{(n)}| t^{n/2}(1+t^{-1/2})^{n-1}\leq Ct^{1/2}
$$
and
\begin{eqnarray*}
\lefteqn{\sup_{|x|\geq t^{-1/2}}|(e^{-tx^2})^{(n)}|(1+|x|)^{n-1}\leq}\\
&\leq &\sup_{|x|\geq t^{-1/2}}(1+|x|)^{n-1}/|x|^{n-1}
\sup_{|x|\geq t^{-1/2}}|(e^{-tx^2})^{(n)}||x|^{n-1}\\
&=&t^{1/2}(1+t^{1/2})^{n-1}\sup_{|x|\geq 1}|(e^{-x^2})^{(n)}||x|^{n-1}\leq C't^{1/2}
\end{eqnarray*} 

where $C$ and $C'$ are constants independent on $t$.
$\;\;\sqcap \!\!\!\! \sqcup$ 
\end{proof}

We now give a combinatorial description of the functions $g_n$.

Denote by $S_n$ the set of subsets $K$ of
$\{0,1,\ldots,n\}$ with the property that $0,n\in K$ and 
$K^c=\{0,1,\ldots,n\}\setminus K$ containes no consecutive integers.
Also denote by 
$\tilde{C}_n=\{(K,\sigma),\, K\in S_n, \sigma \in \{0,1\}^{n}\}$.
The string $\sigma=
\sigma_0\sigma_1\ldots\sigma_{n-1}$, $\sigma_i\in \{0,1\}$ is going to
determine some choices.

To any $(K,\sigma)\in \tilde{C}_n$ 
we are going to associate the continuous
function
$f_{K,\sigma}:[0,\infty)^{n+1}\to 
\Sigma^{(3)}\times \ldots \times\Sigma^{(3)}$ ($n+1$ copies)
defined by 
\begin{equation}
\label{eq-tensor}
f_{K,\sigma}(t_0,\ldots,t_n)=
(e^{-t_0x^2}x^{i_0},e^{-t_1x^2}x^{i_1},\ldots,e^{-t_nx^2}x^{i_n})
\end{equation}
The powers $i_p\in \{0,1,2\}$ are determined
by 
$i_p=j_p+k_p$ where 
\begin{eqnarray}
j_p=\left \{ \begin{array}{ll}
1 & \mbox{ if } p-1,p\in K \mbox{ and } \sigma_{p-1}=1\\
0 & \mbox{ otherwise }
\end{array} \right.
\end{eqnarray}
\begin{eqnarray}
k_p=\left \{ \begin{array}{ll}
1 & \mbox{ if } p,p+1\in K \mbox{ and } \sigma_{p}=0\\
0 & \mbox{ otherwise }
\end{array} \right.
\end{eqnarray}
We agree that $x^0=1$. 
Denote by $L(K)\subset K$ the set of elements $p\in K$ such that
$p+1$ is also in $K$. 
If we denote $m(K)=|L(K)|$ then 
\begin{equation}
\label{eq18}
m(K)=2|K|-n-2=i_0+i_1+\ldots+i_n
\end{equation}
which is proved by induction. 
Here and in what follows $|A|$ denotes the
number of elements of a finite set $A$. Note that the continuity 
follows from Lemma \ref{l-cont}. 

Here is an example. Let $n=7$, $K=\{0,1,3,4,5,6,7\}$ and 
$\sigma=0\sigma_1\sigma_21101$.
Then $L(K)=\{0,3,4,5,6\}$, $i_1=i_2=i_3=i_6=0$, $i_0=i_4=i_7=1$
and $i_5=2$:
$$f_{K,\sigma}(t_0,\ldots,t_7)=(e^{-t_0x^2}\!x,e^{-t_1x^2}\!,e^{-t_2x^2}\!,
e^{-t_3x^2}\!,e^{-t_4x^2}\!x,e^{-t_5x^2}\!x^2
,e^{-t_6x^2}\!,e^{-t_7x^2}\!x)$$
Note that the above function does not depend upon $\sigma_1$ and
$\sigma_2$.

We define  $C_n$ to be the set of equivalence classes in $\tilde{C}_n$
defined by the equivalence relation $\simeq$, where
$(K,\sigma)\simeq (K,\sigma')$ if and only if $\sigma_p=\sigma_p'$ for any 
$p\in L(K)$. 
Let $K\in S_n$, $K\not=\{0\}$, 
we will denote by $K^p=(K\setminus \{0\})-1\in S_{n-1}$ 
if $1\in K$, and $K^p=(K\setminus \{0\})-2 \in S_{n-2}$ if
$1\not\in K$.
Then we define by induction $\nu_{\{0\}}=0$,
$\nu_K=\nu_{K^p}+1$ if $1\not\in K$,
and $\nu_K=\nu_{K^p}+n+1$ if $1\in K$. Let
$$|(K,\sigma)|=\sum_{p\in L(K)}\sigma_p+\nu_K
$$ 
which obviously depends
only on the equivalence class of $(K,\sigma)$. 

We are going to make the convention that,
unless otherwise specified, {\em all integrals are
integrals of continuous functions with values in
the complete locally convex space} $\Om(C^\infty(\Bbb R))$. Also we make
the following conventions.
The sets denoted by 
$K$ give the position of the exponentials, those denoted 
by $L$ give which positions are ocupied by $1$'s and
the sets denoted by $M$ will identify which positions are 
either $e$'s or $e^{-sx^2}x^2$.

Define, using the above notations,
$\eta_n(s;t)\in \Omega^*(C^\infty(\Bbb R))$ by the equation  
\begin{equation} 
\label{eq8}
\eta_n(s;t)= (-1)^n\int_{t\Delta_n}  g_n(s,t_0,\ldots,t_n) dV(t_0,\ldots,t_n)
\end{equation}
and 
\begin{equation}
\label{eq13}
\chi_{K,\sigma}(t)=(-1)^{|(K,\sigma)|}\int_{t\Delta_K}f_{K,\sigma}dV
\end{equation}
where $(K,\sigma)\in C_n$,
$\Delta_K=\{t\in \Delta_n, \supp(t)\subset K\}$.
A consequence of the way 
$f_{K,\sigma}$ was defined is that it
depends only on $\sigma_p$ with $p\in L(K)$, so
by abuse of notation 
we will write from now on $f_{K,\sigma}$ where $(K,\sigma)\in C_n$
is an equivalence class.

Note that from the definitions the largest element
$\max(K)=n$ if $K\in S_n$. 

In the statement of the following  Lemma we are going to use the function
$\Phi:C^\infty(\Bbb R)^{\otimes n+1} \to \Omega^*(C^\infty(\Bbb R))$
defined in Lemma \ref{lemma1}.
\begin{lemma} 
\label{lemma4.7}
Let $\psi_n(t)=\displaystyle{\sum_{(K,\sigma)
\in C_n}}\chi_{K,\sigma}(t)$ and 
$\omega_n(t)=\Phi(\psi_n(t))$. 

(i) The functions $\omega_n(t)$ 
satisfy  $\omega_0(t)=e^{-tx^2}$
and fit into the following reccurence relation
\begin{equation}
\label{eq12}
\omega_{n+1}(t)
=-\int_0^t e^{-sx^2}(dx\omega_{n}(t-s)+dede\omega_{n-1}(t-s))ds,\;\;n\geq 0
\end{equation}
if we let $\omega_{-1}(t)=0$.

(ii) $\sum_{n=0}^\infty \eta_n(s;t)=\sum_{n=0}^\infty s^n\omega_n(t)$.

(iii) $\eta_n(s;t)=\displaystyle{\sum_{|K|=n+1}}
s^{\max(K)}\Phi(\chi_{K,\sigma}(t))$

(iv) The functions $\omega_n(t)$ also satisfy the recurence relation
\begin{equation}
\omega_{n+1}(t)
=-\int_0^t (\omega_{n}(t-s)dx+\omega_{n-1}(t-s)dede)e^{-sx^2}ds
\end{equation}
\end{lemma}
\begin{proof} The proof of (i) is obtained by 
induction as follows.

We get from the definition that 
$e\omega_n(t)e=\omega_n(t)$ for any $n\geq 0$ and $t\geq0$.
Using these relations, the fact that $x$ is
odd and $e(de)e=0$ we see that
the right hand side of the equation (\ref{eq12})
can be written as
$$
\int_0^t e^{-sx^2}(dx\omega_{n}(t-s)+dede\omega_{n-1}(t-s))ds=
$$
$$
=\int_0^t e^{-sx^2}(xd\omega_n(t-s)+
d(x\omega_{n}(t-s))+ded\omega_{n-1}(t-s))ds
$$
Using the induction hypothesis for $\omega_n(t)$ and $\omega_{n-1}(t)$
we obtain the following relations
$$
\int_0^t e^{-sx^2}xd\omega_n(t-s)ds=
-\sum_{(K,\sigma)\in C_n}\Phi(\chi_{K',0\sigma})
$$
$$
\int_0^t e^{-sx^2}d(x\omega_{n}(t-s))ds=
-\sum_{(K,\sigma)\in C_n}\Phi(\chi_{K',1\sigma})
$$
where $K'=\{0\}\cup (K+1)$ 
$$
\int_0^t e^{-sx^2}ded\omega_{n-1}(t-s))ds=
-\sum_{(K,\sigma)\in C_{n-1}}\Phi(\chi_{K'',\sigma''})
$$
where $K''=\{0\}\cup(K+2)$ and $\sigma_{p+2}''=\sigma_p$.

Note that in both cases $(K')^p=K$ and $(K'')^p=K$ and that
$L(K')=\{0\}\cup (L(K)+1)$ and $L(K'')=L(K)+2$.
The above relations justify the choice of the
signs $(-1)^{|(K,\sigma)|}$ in agreement with Lemma \ref{lemma1}.

Since, using the above notation,
$C_{n+1}$ is the disjoint union of the 
sets $\{(K',0\sigma),(K,\sigma)\in C_n\}$,
$\{(K',1\sigma),(K,\sigma)\in C_n\}$ and 
$\{(K'',\sigma''),(K,\sigma)\in C_{n-1}\}$
the above relations immediately give (i).

Write $\sum_{n=0}^\infty \eta_n(s;t)=\sum_{n=0}^\infty \omega_n'(s;t)$
where $\omega_n'(s;t)$ is homogeneous of degree $n$.
The equations (\ref{eq10}) and (\ref{eq8}) show that
$s^n\omega_n(t)$ and $\omega_n'(s;t)$ satisfy the same reccurence relation
with the same initial terms $0$ and $e^{-tx^2}$ so they are
equal. This takes care of (ii). The relation (iii)
is just a restatement of (ii).

The last reccurence relation holds true because it does so for
$\omega_n'(1;t)$ as seen from the equation (\ref{eq9}).\sqep
\end{proof}
\begin{lemma} We have
$\Phi(\chi_{K,\sigma}(t)) \in {\cal F}^{(n)}$
for any $(K,\sigma)\in C_n$.
\end{lemma}
\begin{proof}
Observe that for any fixed 
$k$ and $\epsilon>0$ the the restriction of the 
function $\chi_{K,\sigma}(t_0,t_1,\ldots,t_n)$ to
$t_k\geq \epsilon$
is continuous as a function to ${\cal F}_k^{(n)}$
(${\cal F}^{(n)}_k$ is as defined in the previous 
section, equation (\ref{missing-eq})).\sqep
\end{proof}
From the above two lemmata  we obtain
\begin{corollary}We have
${\eta}_n(s;t),\omega_n(t),\sum_{n\geq0}^\infty\omega_n(t)
\in \overline{\Omega}$ $\forall n \geq 0$, $\forall t>0$.\sqep
\end{corollary}
\begin{definition} We define the exponential
of the (super-)curvature in $\overline{\Omega}$
by the equation
$$
\exp(-t(x+s\nabla)^2)=\sum_{n=0}^\infty s^n\omega_n(t)
\in \overline{\Omega}
$$
\end{definition}
The following Proposition justifies  the previous definition.
\begin{proposition} 
\label{prop3.8}
(i) $\omega_0(t)=e^{-tx^2}$ for any 
$t\geq 0$ and $\omega_n(0)=0$ for $n>0$.

(ii) $\frac{\partial}{\partial t}
\sum_{n=0}^\infty s^n\omega_n(t)=
-(x^2+se(dx)e+s^2edede)\sum_{n=0}^\infty s^n\omega_n(t)$ for any $t\geq 0$
where the derivative is in $\Om(C^\infty(\Bbb R))$.

(iii) $[x,\omega_n(t)]+e(d\omega_{n-1}(t))e=0$.
\end{proposition}
\begin{proof} The first part is a direct consequence of the definitions.

The relation in (ii) is equivalent to 
$$
\frac{\partial}{\partial t}\omega_{n+1}(t)=-x^2\omega_{n+1}(t)
-e(dx)\omega_n(t)-edede\omega_{n-1}(t)
$$ 
which follows by 
differentiating the equations in Lemma \ref{lemma4.7} (i). 
(Use the substitution $s\to t-s$.)

We proceed by induction on $n$ using the result in (ii). 

Let $F(t)=[x,\omega_n(t)]+e(d\omega_{n-1}(t))e=0$.
We  have using the induction hypothesis 
that $\frac{\partial}{\partial t}F(t)=-x^2F(t)$. This shows that
the function $t\to e^{tx^2}F(t)$ has vanishing first derivative
so it is constant. The constant is $0$ as follows from (i).\sqep
\end{proof}

The simplex $\Delta_{\{-1\}\cup K}$ is defined by
$\Delta_{\{-1\}\cup K}=\{(t_{-1},t_0,\ldots,t_n),t_i\geq 0 , \sum t_i=1,
\supp(t)\subset\{-1\}\cup K\}$. Let  
$$\tilde{\chi}_{K,\sigma}(t)=(-1)^{|(K,\sigma)|}\int_{t\Delta_{\{-1\}\cup K}}
e^{-t_{-1}x^2}\otimes h_{K,\sigma}dV$$
\begin{theorem} 
\label{thm4.12}
The functions $\psi_n(t)$  define cycles in Connes'
cyclic complex. More precisely we have the following

(i) $b(\psi_{n+1}(t))=-\tilde{\psi}_{n}(t)$ where 
$\tilde{\psi}_n(t)=\displaystyle{\sum_{(K,\sigma)
\in C_{n-1}}}\tilde{\chi}_{K,\sigma}(t)$.

(ii) $\tilde{\psi}_{n}(t)\in (1-\mbox{\bf t}_{n}) {\cal F}^{(n)}$
where $\mbox{\bf t}_{n}(a_0\otimes \ldots \otimes a_{n})=(-1)^{n+\nu}
a_{n}\otimes a_0 \otimes \ldots \otimes a_{n-1}$ defines the action of 
the cyclic group $\Bbb Z_{n+1}$ on ${\cal F}^{n}$ for 
$\nu=\partial a_{n}(\partial a_0+\ldots +\partial a_{n-1})$.
\end{theorem}

We are going to use Lemma \ref{lemma1} quite often so
some comments are in order. 

According to the  definition
$\Omega^*(C^\infty(\Bbb R))$ containes two copies of 
the projective tensor product 
$C^\infty(\Bbb R)^{\otimes n+1}$, one of them is 
the image of $\Phi$ and the other one is the image of 
$d\circ\Phi$. It is important to point out that 
$\Phi$ is not the identity but introduces some signs
according to the definition in Lemma \ref{lemma1}.
Also we are going to use the same symbol $\Phi$ for different
values of $n$. This should cause no confusion.

The grading automorphism $\alpha$ on each component 
of the tensor products gives
rise, by means of $\Phi$, to a grading $P$ on 
$\Omega^*(C^\infty(\Bbb R))$
which is different however from the standard grading 
denoted also $\alpha$ exactly
because $d$ is odd. The precise relations are given by
$P(\xi)=(-1)^n\alpha(\xi)=(-1)^{n+\partial \xi}\xi$
for $\xi\in \Omega^n(C^\infty(\Bbb R))$. 
Thus $dx$ is even for the $\alpha$ grading but 
odd for $P$. Also $\alpha(\omega_n(t))=\omega_n(t)$
but $P(\omega_n(t))=(-1)^n\omega_n(t)$.

The reason for introducing $P$ is the following equation
\begin{equation}
\label{eq-sign}  
\Phi^{-1}(adb)=P(\Phi^{-1}(a))\otimes b\;\;\forall a\in \mbox{Im}(\Phi)\subset
\Omega^*(C^\infty(\Bbb R)), b\in C^\infty(\Bbb R)
\end{equation}

We are ready now to start the proof.
\paragraph{Proof of the Theorem.}
Let $\Omega_{n+1}(t)=\omega_n(t)dx+\omega_{n-1}(t)dede$.
Lemma \ref{lemma4.7} (iv) gives that
$$
-\Phi(\psi_{n+1}(t))=\int_0^te^{-sx^2}\Omega_{n+1}(t-s)ds+
\int_0^t[\Omega_{n+1}(t-s),e^{-sx^2}]ds
$$
Using then Lemma \ref{lemma1} we obtain that 
$$
-\psi_{n+1}(t)=\int_0^te^{-sx^2}\Phi^{-1}(\Omega_{n+1}(t-s))ds
+b\int_0^t\Phi^{-1}(\Omega_{n+1}(t-s))\otimes e^{-sx^2}ds
$$
Using again Lemma \ref{lemma1}, the definition 
of $\Omega_{n+1}(t)$ and equation (\ref{eq-sign})
we obtain 
$$
-\Phi(b\psi_{n+1}(t))=(-1)^n\int_0^te^{-sx^2}([P(\omega_n(t-s)),x]
+[P(\omega_{n-1}(t-s))de,e])ds
$$
which finally gives using Proposition \ref{prop3.8} (iii)
\begin{eqnarray*}
\lefteqn{(-1)^{n+1}\Phi(b\psi_{n+1}(t))}\\
&=&-P\int_0^te^{-sx^2}([\omega_n(t-s),x]-\omega_{n-1}(t-s)de)ds\\
&=&-P\int_0^te^{-sx^2}d\omega_{n-1}(t-s)ds\\
&=&-P\Phi(\tilde{\psi}_{n}(t))=(-1)^n\Phi(\tilde{\psi}_n(t))
\end{eqnarray*}
This proves (i).
In order to prove (ii) we begin by making some remarks.

Denote for $(K,\sigma)\in C_n$ by $\phi(K,\sigma)=\{p\in\{0,1,\ldots,
n\}: i_p=1\}\subset K$ where $i_p$ are the exponents appearing in the definition
of $f_{K,\sigma}$, see equation (\ref{eq-tensor}).

{\em Claim 1.} $n-|\phi(K,\sigma)|$ is even. Given any subset
$L\subset K$ such that $n-|L|$ is even there exists at most one
equivalence class $(K,\sigma)\in C_n$ such that $\phi(K,\sigma)=L$.
Order increasingly $\{0,1,\ldots,n\}\setminus L=\{p_0,p_1,\ldots,p_{2l}\}$,
then a solution exists if and only if $M(L)=\{p_0,p_2,p_4,\ldots,p_{2l}\}
\subset K$.

{\em Proof of the Claim.}
We know that $i_0+i_1+\ldots+i_n=2|K|-n-2$ (equation \ref{eq18}) 
from which it immediately
follows that the number of $1$'s has the same parity as $n$.

Suppose first  that $K=\{0,1,\ldots,n\}$. We have the relation
\begin{equation}
\label{eq-det1}
i_0+i_1+\ldots+i_q=q+1-\sigma_q\;,\,q=0,1,\ldots,q-1
\end{equation}
The above equation shows that $\sigma_q=\sigma_{q-1}$ if
$i_q=1$, that $\sigma_{q-1}=0$ and $\sigma_q=1$ if $i_q=0$,
and that $\sigma_{q-1}=1$ and $\sigma_q=0$ if $i_q=2$. We also
see that $K\setminus L$ tells us where changes in the string
$\sigma$ occur. From this it follows that, after ignoring the 
exponents equal to $1$,
the positions of the exponents equal to $0$ alternate with
those equal to $2$. The above equation tells us that 
since $1+1+\ldots+1+2$ is ``too big'', we necesssarily
have to begin the alternation procedure with a $0$. 
This determines all $i_q$'s and hence all $\sigma_q$'s again
in view of the previous equation.
This proves the existence and
uniqueness in the case $K=\{0,1,\ldots,n\}$.

In order to determine uniqueness and when solutions exist
for an arbitrary given $K$ 
observe that if $(K_1,\sigma)\in C_n$ and 
$\sigma_{p-1}=1$ and $\sigma_p=0$ then 
$K=K_1\setminus \{p\}$ is still in $S_n$ and 
$\phi(K,\sigma)=\phi(K_1,\sigma)$. Thus we can
add elements to $K$ if we appropriately define the
previously undefined $\sigma$'s. 

This proves
the uniqueness. Moreover to every element that we 
add to $K$ we get the  corresponding exponent to change from 
$0$ to $2$ and all the other ones to stay put. 
This means that in order for the above procedure to 
allow us to decrease the number of elements in $K$ we 
can take out only elements that have the corresponding
exponent equal to $2$. If we examine the existence in
the case $K=\{0,1,\ldots,n\}$ we see that the above
positions are $p_1,p_3,\ldots,p_{2l-1}$.

This proves the claim.

Note for further reference that if $(K,\sigma)$ and $(K_1,\sigma)$
are as above then $\nu_{K}=\nu_{K_1} \pmod 2$ and hence
\begin{equation}
\label{eq-sign2}  
|(K,\sigma)|=|(K_1,\sigma)|+1 \pmod 2
\end{equation}

Consider three disjoint sets $L,M,M_1\subset \{-1,0,1,\ldots,n\}$  such
that $n-|L|$ is even and $|M|=(n-|K|)/2$.

We define 
$$
\zeta_{L,M,M_1}(t_{-1},t_0,\ldots,t_n)=f_{-1}(t_{-1})\otimes
f_0(t_{0})\otimes \ldots \otimes f_n(t_n)
$$
where
\begin{eqnarray*}
f_p(s)=\left \{ \begin{array}{ll}
e^{-sx^2} & \mbox{ if } p\in M\\
e^{-sx^2}x & \mbox{ if } p\in L\\
e  & \mbox{ if } p \in M_1\\
e^{-sx^2}x^2 & \mbox{ otherwise }
\end{array} \right.
\end{eqnarray*}
Also consider 
$$
\tilde{\zeta}_{L,M,M_1}(t)=\int_{t\Delta_{S}}
\zeta_{L,M,M_1}dV
$$
where $S=\{-1,0,1,\ldots,n\}\setminus M_1$.

{\em Claim 2.} $\zeta_{L,M}=
\sum_{M_1\subset M}(-1)^{|M_1|} \tilde{\zeta}_{L,M,M_1}$ 
where the sum is over all subsets of $\{-1,0,1,\ldots,n\}\setminus (L\cup M)$,
does not depend on $M$.

{\em Proof of the Claim.} Since all the sets $M$ have the
same number of elements in $\{0,1,\ldots,n\}\setminus L$ it is
enough to show that $\zeta_{L,M}(t)=\zeta_{L,M'}(t)$ for
$M$ and $M'$ differing by exactly one element. More
precisely we assume that $M=M''\cup \{j\}$ and 
$M'=M''\cup \{k\}$ (the case when $M$ must be empty is trivial). 

The difference 
$\zeta_{L,M}(t)-\zeta_{L,M'}(t)$ can then be written
as the sum of all terms of the form 
$$
(-1)^{|M_1|}(-\tilde{\zeta}_{L,M,M_1\cup \{k\}}
-\tilde{\zeta}_{L,M',M_1\cup \{j\}}
+\tilde{\zeta}_{L,M,M_1}
+\tilde{\zeta}_{L,M',M_1})
$$
where $M_1$ ranges through all subsets of 
$\{-1,0,1,\ldots,n\}\setminus (L\cup M''\cup \{k,j\})$.
A simple integration with respect to the
extra variable shows that the above term is $0$.
(Use the same argument as in \cite{NiD}, equation (2.5),
page 450.) 

The proof of the Claim is now complete.

The above claim tells us that we can define
$\zeta_L(t)=(-1)^{\Sigma(L)}\zeta_{L,M}(t)$ 
where $\Sigma(L)=\sum_{k\in L}(n-k)$, for an arbitrary choice 
of $M$.

{\em Claim 3.} $\tilde{\psi}_{n}(t)=(-1)^n\sum \zeta_L(t)$
where $L$ runs through all subsets of $\{0,1,\ldots,n\}$
with the property that $n-|L|$ is even.

{\em Proof of the Claim.} 
It follows from the first claim that 
$$
e^{t_1}\otimes f_{K,\sigma}(t_0,\ldots,t_n)=
\zeta_{L,M(L)\cup\{-1\},M_1}(t_{-1},t_0,\ldots,t_n)
$$ 
for $L=\phi(K,\sigma)$ and $M_1=\{0,1,\ldots,n\}\setminus K$.
Moreover $|(K,\sigma)|=n+\Sigma(L)+|M_1| \pmod 2$, this follows from
equation (\ref{eq-sign2}) by induction on $|M_1|$. For $M_1=\emptyset$
we use $\nu_K=(n+1)(n+2)/2-1$ if $K=\{0,1,\ldots,n\}$ and
$\sigma_0+\sigma_1+\ldots+\sigma_{n-1}=\sum{k=0}^n (n-k)i_k$. 
Since the last sum coincides with $\Sigma(L) \pmod 2$ this shows that
in the next two sums all the terms coincide:
$$
\sum_{(K,\sigma),\phi(K,\sigma)=L}
\tilde{\chi}_{K,\sigma}(t)=\zeta_{L,M(L)\cup \{-1\}}(t)
$$
From this the proof of the claim follows.

In the following Claim  $\mbox{\bf t}_{n+1}$ also denotes
the cyclic permutation of the set $\{-1,0,1,\ldots,n\}$:
$\mbox{\bf t}_{n+1}(n)=-1$, $\mbox{\bf t}_{n+1}(-1)=0$,
$\mbox{\bf t}_{n+1}(0)=1$ ,..., $\mbox{\bf t}_{n+1}(n-1)=n$

{\em Claim 4.} Suppose $L,L_1 \subset \{0,1,\ldots,n\}$
are such that $L=\mbox{\bf t}_{n=1}^{k+1} L_1$ where $k>0$ is least with this 
property. Then $\zeta_{L}=-\mbox{\bf t}_{n+1}^{k+1}\zeta_{L_1}$.

{\em Proof of the Claim.} We observe that if $L'$ is an arbitrary 
subset of $\{-1,0,\ldots,n\}$ such that $|L|=n \pmod 2$ then if 
we denote $L''=\mbox{\bf t}_{n+1}^{-1}$
we have
$\mbox{\bf t}_{n+1}\zeta_{L''}=\zeta_L$ for $-1 \in L$ and
$\mbox{\bf t}_{n+1}\zeta_{L''}=\zeta_L$ for $-1 \not\in L$.
This is seen in the following way. In the first case $\Sigma(L'')=
\Sigma(L')\pmod 2$ and there is no sign induced by the cyclic permutation
(because an odd element ``jumps'' over other $|L|-1$ odd elements).
In the second case we have $\Sigma(L'')=\Sigma(L)+|L|$ and the sign
of the cyclic permutation is the same as in the trivialy 
graded case (as no odd element ``jumps''), that is $(-1)^{n+1}$.

This proves the last Claim as well as the Theorem.\sqep

\paragraph{ }
We proceed now to state and prove our main theorem after reviewing
the main concepts and assumptions involved in its statement.

\paragraph{Assumptions and notation.} We assume that we are given the following:
\begin{enumerate}
\item a cycle $\cycle$ (definition \ref{def1}), 
that is $\nabla:\Omega \to \Omega$ is a 
graded derivation such that $\nabla^2$ is inner and the
trace $\tau$ vanishes on covariant derivatives $\tau(\nabla(a))=0$.
\item the larger algebra 
$\tilde{\Omega}=\Omega\oplus \Omega X\oplus X\Omega \oplus X\Omega X$,
as in equation (\ref{eqOmega}), together with the canonical differential $d$ and
extension $\tilde{\tau}$ of the trace $\tau$ such that $\cycletil$
is an exact cycle.
\item we assume $\tau(F_{-2m-1}\Omega)=0$ and we will use the $2m$-cyclic
cocycle $\varphi^\tau$ on $A=\Omega/F_{-1}\Omega$ corresponding 
to $\cycle$ and $2m$. See definition \ref{def2.7.}.

\item a generalized $\theta$-summable Fredholm operator 
in $A$, that is a continuous morphism $D:{\cal D}\to A$,
with `index'  $\mbox{Ind}(D) \in K_0(A\rtimes \Bbb Z_2)$ defined by the equation
(\ref{eqIND}).

\item we also assume that $\cycle$ is finite summable (over ${\cal D}$), 
definition \ref{defFS}, that is that the morphism $D$ extends to a morphism
$\overline{D}:\overline{\Omega} \to \tilde{\Omega}$ and that
$\tilde{\tau}$ vanishes on $\overline{D}([\Sigma^{(\infty)},\overline{\Omega}])$.
\end{enumerate}

The cocycle $\varphi^\tau$ provides us then
with an  $2m$-cyclic cocycle $\psi^\tau$ on $A\rtimes \Bbb Z_2$ (equation
(\ref{eq6-})), and hence to a
morphism, the Connes-Karoubi character, \cite{Co2,Karoubi2}
\begin{equation}
\psi^\tau_*:K_0(A\rtimes \Bbb Z_2) \to \Bbb C
\end{equation}
Our main problem is to compute the Connes-Karoubi character of $D$, more
precisely $\psi_*^\tau(\mbox{Ind}(D))$. The main theorem
will give an expression for this `characteristic number' in case
$\psi$ satisfies certain further continuity conditions (finite summability).

Here $\Sigma^{(\infty)}=\cup \Sigma^{(n)}$ is an algebra of polynomial-like functions
defined in equation (\ref{eq10}), and $\overline{\Omega}$ is a certain
completion of $\Omega^*({\cal D})$, the universal differential graded algebra
of ${\cal D}$, see equation (\ref{eqcalD}).

The exponential of the supercurvature is defined
by 
$$
exp(-t(D+s\nabla)^2)=\overline{D}(\sum_{n=0}^\infty s^n \omega_n(t))
$$

\begin{theorem}
\label{thmMain} Using the above notation and assumptions for the cycle
$\cycle$ and the generalized $\theta$-summable operator 
$D:{\cal D} \to A=\Omega/F_{-1}\Omega$ we have that
\begin{equation}
\tau(exp(-t(D+s\nabla)^2))=\!P(st)\!
=(st)^{2m}\psi_*^\tau(\mbox{Ind(D)})+\mbox{`lower order 
terms'}
\end{equation}
where $P$ is a polynomial of degree $2m$ with constant coefficients and
$\psi^\tau$ is the cyclic cocycle associated to $\tau$ and $2m$, assuming
that $\tau$ vanishes on $F_{-2m-1}\Omega$.

All the lower terms of the polynomial $P$ vanish if $\Omega$ is graded
and $\nabla$ is of degree $1$.
\end{theorem}
\begin{proof} We are going to use the notations explained above.

The linear map $\tilde{\tau}\circ \overline{D}$ defines a closed graded
trace on $\overline{\Omega}$. Let $\Phi$ be the map defined in Lemma 
\ref{lemma1}, with components $\Phi_n(a_0\otimes a_1 \ldots a_n)=
(-1)^\mu a_0 da_1 \ldots da_n$ then 
$\xi_n=\tilde{\tau}\circ \overline{D}\circ 
\Phi_n$ are cyclc cocycles on ${\cal D}$
for all $n$ and $\varphi^\tau=\xi_{2m}$ is the cyclic cocycle in the statement
of the theorem. 

This gives $\omega_n(t)=\Phi_n(\psi_n(t))$ using the notation 
of lemma \ref{lemma4.7}, and hence
$$
\tau(exp(-t(D+s\nabla)^2))=\sum_{n=0}^m s^{2n}\xi_{2n}(\psi_{2n}(t))
$$
By theorem \ref{thm4.12} we know that $\psi_{2n}(t)$ are cyclic cyles.
The components $ch_{2n}(\mbox{Ind}(id)) \in 
HC_{2n}({\cal D}\rtimes \Bbb Z_2)_{odd}
\simeq HC_{2n}({\cal D})_{even}$ 
of the Connes-Karoubi character are also cyclic cycles. Since the
homology of $({\cal F})_{n\geq 0}$ is singly generated, Theorem \ref{thm3.6},
it follows that $\psi_{2n}(t)=f_{2n}(t)ch_{2n}(\mbox{Ind}(D))$
in homology, for some complex functions $f_{2n}$. This gives 
$$
\xi_{2n}(\psi_{2n}(t))=f_{2n}(t)\xi_{2n}(ch_{2n}(\mbox{Ind}(id)))
$$
We see then that the theorem is equivalent to 
$f_{2n}(t)=t^{2n}$. Since these functions do not depend on 
any choice, we can check our computation in the simplest example
when the operator $D$ vanishes, more precisely, when 
$D(f)=f(0)$. We can also assume that $\Omega=\Omega^*(M,E)$, the
space of smooth forms with coefficients in the smooth bundle
$E$ on a  smooth compact manifold $M$ of dimension $2n$.
In this case the computation reduces to the 
Chern-Weyl definition of Chern classes in terms of 
connection and curvature.\sqep
\end{proof}

The following example will spell in more detail the relation
between our theorem and the Chern-Weyl construction of 
characteristic classes.
We will use Proposition \ref{prop4.4}.
We make the observation that the extension 
$\overline{\rho}:\Sigma^{(\infty)} \to {\frak B}$
need not be unital and hence that the value of 
$exp(-t(D+s\nabla)^2)$ at $t=0$ is not the unit of ${\frak B}$
but rather an idempotent $\overline{e}$, the strong
limit of $exp(-t(D+s\nabla)^2)$ as $t$ decreases to $0$.

Consider a smooth compact manifold $X$ and $p_0,p_1$ two 
projections in $M_n(C^\infty(X))$ and define 
$$
\overline{\rho}(f)=f(0)p_0\oplus f(0)p_1 \in A=M_{2n}(C^\infty(X))
$$
where $A$ is endowed with the grading defined by the
matrix $\gamma=1\oplus (-1)$. Define ${\frak B}=\Omega=M_{2n}(\Omega^*(X))$.

Then we have $A\rtimes \Bbb Z_2 \simeq A\oplus A$ where the first
map sends $v$ to $\gamma$ and the second one
sends $v$ to $-\gamma$. The index is given by the odd element
$$
\mbox{Ind}(D)=([p_0]-[p_1],-[p_0]+[p_1]) \in K_0(A) \oplus K_0(A).
$$

The even super traces on the algebra $\Omega$ correspond to currents on $X$:
$\tau(\omega)=\langle Str(\omega),c\rangle $ where $c \in \Omega^{2m}(X)'$,
the dual space. The cocycles $\varphi^\tau$ and $\psi^\tau$ will
be given by the formulae
$$
\varphi^\tau (a_0,a_1,\ldots,a_{2m})
=\langle Str(a_0da_1\ldots da_{2m}),c\rangle  \, , \;\; a_0,\ldots a_{2m} \in A
$$
and
$$
2\psi^\tau (b_0,b_1,\ldots,b_{2m})
=\langle (tr\oplus -tr)(b_0db_1\ldots db_{2m}),c\rangle  \, , \;\; b_0,\ldots b_{2m} \in A \oplus A
$$
which give:
\begin{eqnarray*}
\psi^\tau(\mbox{Ind(D)}) & = & \langle tr(p_0dp_0dp_0)^{m} - tr(p_1dp_1dp_1)^{m},c\rangle /(m!)\\
& = & (2\pi \imath)^m \langle  ch(p_0)-ch(p_1),c\rangle .
\end{eqnarray*}
Since 
$exp(-t(D+s\nabla)^2)=exp(s^2t^2 p_0dp_0dp_0)p_0 \oplus exp(s^2t^2 p_1dp_1dp_1)p_1$
we obtain 
$$
\tau (exp(-t(D+s\nabla)^2)) = 
(st)^{2m} \langle Str (p_0dp_0dp_0)^{m} \oplus (p_1dp_1dp_1)^{m},c\rangle /(m!)
$$
and this checks with the theorem.


\end{document}